\documentclass[fleqn]{zzz}
\usepackage{mathrsfs}
\usepackage{tikz}
\usepackage[pdftex]{hyperref}

\hypersetup{
	colorlinks = true,
	urlcolor = blue,
	linkcolor = red!60!black,
	citecolor = green!60!black
}

\newtheorem{thm}{Theorem}[section]
\newtheorem{cor}[thm]{Corollary}

\newtheorem{rem}[thm]{Remark}
\newtheorem{lem}[thm]{lem}
\newtheorem{defn}[thm]{Definition}
\newtheorem{prop}[thm]{Proposition}

\makeatletter
\def\eqnarray{\stepcounter{equation}\let\@currentlabel=\theequation
	\global\@eqnswtrue
	\tabskip\@centering\let\\=\@eqncr
	$$\halign to \displaywidth\bgroup\hfil\global\@eqcnt\z@
	$\displaystyle\tabskip\z@{##}$&\global\@eqcnt\@ne
	\hfil$\displaystyle{{}##{}}$\hfil
	&\global\@eqcnt\tw@ $\displaystyle{##}$\hfil
	\tabskip\@centering&\llap{##}\tabskip\z@\cr}

\def\endeqnarray{\@@eqncr\egroup
	\global\advance\c@equation\m@ne$$\global\@ignoretrue}

\def\@yeqncr{\@ifnextchar [{\@xeqncr}{\@xeqncr[5pt]}}
\makeatother

\begin{document}

\renewcommand{\PaperNumber}{***}

\FirstPageHeading

\ShortArticleName{On polynomials orthogonal with respect to
an inner product involving higher order differences. The Meixner case}

\ArticleName{On polynomials orthogonal with respect to
an inner product involving higher order differences. The Meixner case}

\Author{Roberto S. Costas-Santos,
Anier Soria-Lorente$^{1}$ and Jean-Marie Vilaire$^{2}$}

\AuthorNameForHeading{Roberto S. Costas-Santos, 
Anier Soria-Lorente, Jean-Marie Vilaire}
\URLaddressD{\href{http://www.rscosan.com}{http://www.rscosan.com}
}
\EmailD{rscosa@gmail.com} 
\Address{$^1$ Facultad de Ciencias T\'ecnicas, Universidad de Granma,
Km. 17.5 de la carretera de Bayano-Manzanillo, Bayamo, Cuba}
\EmailD{asorial@udg.co.cu, asorial1983@gmail.com}

\Address{$^2$ Institut des Sciences, des Technologies et des \'{E}tudes Avanc\'{e}%
es d'Ha\"{\i}ti, \# 10, Rue Mercier-Laham, Delmas 60, Musseau, Port-au-Prince, Haiti.
B.P. 15953}
\EmailD{jeanmarievilaire@yahoo.fr}

\ArticleDates{Received: 15 April 2022 Accepted: 3 June 2022 Published: 6 June 2022}
\Abstract{In this contribution we consider sequences of monic polynomials 
orthogonal with respect to the Sobolev-type inner product 
$$
\left\langle f,g\right\rangle= 
\langle {\bf u}^{\tt M},fg\rangle
+\lambda \mathscr T^j f (\alpha)\mathscr T^{j}g(\alpha),
$$
where ${\bf u}^{\tt M}$ is the Meixner linear operator, 
$\lambda\in\mathbb{R}_{+}$, $j\in\mathbb{N}$, $\alpha \leq 0$, and
$\mathscr T$ is the forward difference operator  $\Delta$, or the backward 
difference operator $\nabla$.
\newline
We derive an explicit representation for these polynomials. 
The ladder operators associated with these polynomials are obtained, and 
the linear difference equation of second order is also given. In addition, 
for these polynomials we derive a $(2j+3)$-term recurrence relation. 
Finally, we find the Mehler-Heine type formula for the particular case $\alpha=0$.
}


\Classification{Primary 33C47; Secondary 39A12}

\section{Introduction}
The Meixner orthogonal polynomials, usually denoted in the literature as 
$M_{n}(x;\beta,c)$, constitute a family of classical orthogonal
polynomials, introduced by J. Meixner in 1934, in his
seminal paper \cite{meixner1934}. When $\beta>0$ and $0<c<1$ they 
are orthogonal with respect to the well-known negative binomial distribution 
of the probability theory, i.e. in such a case
\[
\sum_{x=0}^\infty \binom{x+\beta-1}{x} c^x\, M_{n}(x;\beta,c)M_m(x;\beta,c)=0, \quad m\ne n,\ m, n=0, 1, \dots
\]
so the Meixner linear functional is
\[
 {\bf u}^{\tt M}=\sum_{x=0}^\infty \binom{x+\beta-1}{x} c^x\,\delta_x.
\]
So, they are orthogonal on the uniform lattice in the interval $[0,+\infty )$ and 
they satisfy an hypergeometric-type difference equation on the aforesaid uniform
lattice. 
Because their classical character, their finite differences constitute as well 
an orthogonal polynomial family, their corresponding
orthogonality weight satisfy a Pearson-type difference equation, and even
they satisfy two different kinds structure relations. A distinctive and
interesting characteristic of this family is that they have certain dual
character, that is, every formula one can derives for 
$M_{n}(x;\beta,c)$ has a dual formula with $x$ and $n$ interchanged 
\[
c^{m-n}n!(1+\beta)_{m-1}M_{n}(m;\beta,c)=
c^{n-m}m!(1+\beta)_{n-1}M_m(n;\beta,c),
\]
(see, for example \cite[Ch. VI, Sec. 3]{Chi-78}, \cite[Sec. 6.1]{Ismail-05}, 
\cite[Sec. 2.4] {NikUvSu-91}, and the references therein).

On the other hand, since the first paper \cite{althammer1962} on Sobolev 
orthogonal polynomials published by Althammer, until the present time, 
the results connected to these polynomials have attracted the attention 
of several mathematicians. The name of Sobolev orthogonal polynomials 
was given to those families of polynomials orthogonal with respect to inner 
products involving positive Borel measures supported on infinite subsets 
of the real line, and also involving regular derivatives. Moreover, in the c
ase that the derivatives appear only on function evaluations on a finite 
discrete set, the corresponding families are called \textit{Sobolev-type} 
or \textit{discrete Sobolev} orthogonal polynomial sequences. 
For a recent and comprehensive survey on the subject, see \cite{MX-EM15} 
and the references therein. In the last decade of the past century, 
H. Bavinck introduced the study of orthogonal polynomials with respect 
to the inner product involving differences instead of derivatives 
\begin{equation}\label{IPB}
\langle f,g\rangle_{\lambda}=\int_{\mathbb{R}}f(x)g(x) d\psi (x)
+\lambda(\Delta f)(c)(\Delta g)(c),
\end{equation}
where $\lambda \in\mathbb{R}^{+}$, $c\in\mathbb{R}$ and $\psi $ 
is a distribution function with infinite spectrum, see \cite{B-JCAM95,B-AA95}. 
Moreover, in these works Bavinck obtained algebraic properties and some 
results connected to the location of the zeros of the orthogonal polynomials 
with respect to the inner product \eqref{IPB}. On the other hand, in 
\cite{B-AA95} he proved that the orthogonal polynomials with respect to 
inner product defined in equation\eqref{IPB} satisfy a five term recurrence 
relation. 
Furthermore, in\cite{Bavinck-94} the author considered the inner product 
\begin{equation}
\left\langle f,g\right\rangle =(1-c)^{\beta }\sum_{x=0}^\infty f(x)g(x)
\frac{c^{x}\Gamma (\beta+x)} {\Gamma(\beta)\Gamma(x+1)}
+\lambda f(0)g(0), \label{IPB2}
\end{equation}
where $\beta >0$, $0<c <1$, $\lambda >0$ and $\mathbb{P}$ denote the
linear space of all polynomials with real coefficients. Here, he obtained a
second order difference equation with polynomial coefficients, which the
orthogonal polynomials with respect to \eqref{IPB2} satisfy. 
Then, in\cite{Bavinck-98} the author showed that the Sobolev type Meixner 
polynomials orthogonal with respect to the inner product 
\begin{equation*}
\left\langle f,g\right\rangle =(1-c)^{\beta }\sum_{x=0}^\infty f(x)g(x)\frac{
c^{x}\Gamma (\beta+x)} {\Gamma(\beta)\Gamma(x+1)}+Mf(0)g(0)
+N(\Delta f)(0)(\Delta g)(0),
\end{equation*}
where $\beta >0$, $0<c <1$ and $M, N\geq 0$, are eigenfunctions of a
difference operator. Other results a little most recent, connected with the 
Sobolev Meixner polynomials, can be found in \cite{KD-AAS12,M-JCAM15}.

The structure of the paper is the following: In Section 2, we introduce some
preliminary results about Meixner polynomials which will be very useful in
the analysis presented. In Section 3, we obtain the connection formula
between the Meixner polynomials and the polynomials orthogonal with respect
to the Sobolev-type inner product
\begin{equation}\label{SobIP}
\left\langle f,g\right\rangle= 
\langle {\bf u}^{\tt M},fg\rangle
+\lambda \mathscr T^j f (\alpha)\mathscr T^{j}g(\alpha),
\end{equation}
where ${\bf u}^{\tt M}$ is the Meixner linear operator, 
$\lambda\in\mathbb{R}_{+}$, $j\in\mathbb{N}$, $\alpha \leq 0$, and
$\mathscr T$ is the forward or the backward difference operator,
as well as we deduce the hypergeometric representation of
such polynomials. In Section 4, we find the ladder (creation and
annihilation) operators for the sequence of orthogonal polynomials of
Sobolev type. As a consequence, the second order linear difference equation
associated with them is deduced. And on the other hand, in Section 5, we
determine the $(2j+3)$-term recurrence relation that these polynomials
satisfy. Finally, in Section 6, we determine the Mehler-Heine type formula
for the especial case $\alpha =0$. Indeed, the techniques used in Sections
3, 4 and 5 are based on those used in\cite{Alvarez000,Edy01,Ismail-05},
respectively.
\section{Preliminaries}
We adopt the following set notations: $\mathbb N_0:=\{0\}\cup\mathbb 
N=\{0, 1, 2, 3, \ldots\}$, and we use the sets $\mathbb Z$, $\mathbb R$, 
$\mathbb C$ which represent the integers, real numbers and complex 
numbers respectively. $\mathbb P$ denotes the vector space of univariate, 
complex-valued, polynomials, and let $\mathbb P'$ denote its algebraic dual space. 

We also adopt the following notation and conventions.
We denote by $\langle {\bf u}, p\rangle$ the duality bracket for ${\bf u}\in
\mathbb P'$ and $p\in\mathbb P$.

\begin{defn} For ${\bf u}\in\mathbb P$, $\pi \in\mathbb P$, and $c\in\mathbb C$, 
let $\pi {\bf u}$, $(x-c)^{-1}{\bf u}$, $\nabla {\bf u}$, and $\Delta {\bf u}$ 
be the linear functional defined by
\begin{eqnarray*}
\langle \pi{\bf u}, p\rangle&:=&\langle {\bf u}, \pi p\rangle, \quad p\in\mathbb P,\\
\langle (x-c)^{-1} \bf u, p\rangle&:=&\left\langle {\bf u}, \frac{p(x)-p(c)}{x-c}\right\rangle, 
\quad p\in\mathbb P, \\
\langle \nabla \bf u, p\rangle&:=&-\left\langle {\bf u}, \Delta p\right\rangle,
\end{eqnarray*}
and thus 
$\langle \Delta \bf u, p\rangle:=-\left\langle {\bf u}, \nabla p\right\rangle$,
where $\nabla$ and $\Delta$ are the backward and forward difference operator 
defined as: 
\[
\Delta f(x):=f(x+1)-f(x),\qquad \nabla f(x):=f(x)-f(x-1).
\]
The Dirac delta functional, $\delta_c$, is the functional defined by 
$\langle \delta_c, p\rangle:= p(c)$, $p\in\mathbb P$.
\end{defn}
In order to obtain our derived identities, we rely on properties 
of the Pochhammer symbol (shifted factorial). 
For any $n\in\mathbb N_0$, $a\in\mathbb C$, the Pochhammer symbol is 
defined as
\[
(a)_n:=(a)(a+1)\cdots(a+n-1),\quad n\in\mathbb N_0,
\]
Furthermore, define for all $a,b\in\mathbb C$, 
\[
(a)_b:=\dfrac{\Gamma(a+b)}{\Gamma(a)},
\]
where $a+b\not\in-\mathbb N_0$, and 
we will also use the common notational product convention
\[
(a_1, \dots, a_k)_b:=(a_1)_b\cdots(a_k)_b,
\]
The hypergeometric series, which we will often use, is defined for
$z\in\mathbb C$ such that $|z|<1$, $s,r\in\mathbb N_0$, 
as \cite[(1.4.1)]{Koekoeketal}
\begin{equation}
_{r}F_{s}\left(\begin{array}{c} a_1, \dots, a_r\\ 
b_1, \dots, b_s\end{array};z\right)
:=\sum_{k=0}^\infty \frac{(a_1, \dots, a_r)_k}{(b_1, \dots, b_s)_k}
\dfrac{z^k}{k!}. \label{HSDE}
\end{equation}
where, of course, the parameters must be such that the denominator 
factors in the terms of the series are never zero.
\subsection{The Meixner polynomials}
Let $\beta$ and $c$ be two complex numbers such that $c\ne 0, 1$ and 
$\beta$ is not a negative integer. 
We write $\{M_{n}(x;\beta,c)\}_{n\geq 0}$ for the sequence of Meixner 
polynomials defined by \cite{Chi-78}
\[
M_{n}(x;\beta,c)=\frac{c^n\, n!}{(1-c)^n} \sum_{j=0}^n c^{-j}\binom{x}{j}
\binom{-x-\beta}{n-j}.
\]
These polynomials are orthogonal with respect to the 
linear functional ${\bf u}^{\tt M}\in\mathbb P'$ which is a classical 
functional since it fulfills the Pearson difference equation
\[
\Delta \left(x {\bf u}^{\tt M}\right)=\left(x(c-1)+\beta c\right){\bf u}^{\tt M},
\]
which is equivalent to the Pearson difference equation
\[
\nabla \left(c(x+\beta) {\bf u}^{\tt M}\right)=\left(x(c-1)+\beta c\right){\bf u}^{\tt M}.
\]
\begin{rem}
Observe that when $\beta>0$ and $0<c<1$, then 
\[
\langle {\bf u}^{\tt M},f\rangle:=\sum_{x=0}^\infty \binom{x+\beta-1}{x} c^x\,f(x),
\]
which is a positive definite linear functional. Observe this definition can 
be extended to $|c|<1$ and $\beta\in\mathbb C$ nor a negative integer. 
Moreover, since 
\begin{equation} \label{meorex}
M_{n}(x;\beta,c)=(-1)^n M_{n}(-x-\beta;\beta,c^{-1}),
\end{equation} 
then one can extend the Meixner functional for $|c|>1$.
In\cite[Proposition 9]{CoSa} the authors obtained an integral representation 
for this operator for $\beta,c \in\mathbb C$, with $c\not \in[0,\infty]$ and 
$-\beta\not \in\mathbb N$:
\[
\langle {\bf u}^{\tt M},f \rangle=\int_C
f(z)\Gamma(-z)\Gamma(\beta+z)(-c)^z\,dz,
\]
where $C$ is a complex contour from $-\infty i$ to $\infty i$ separating the increasing 
poles $\{0, 1, 2, \dots\}$ from the decreasing poles $\{-\beta , -\beta -1,-\beta-2,\dots\}$.
\end{rem}
When $c\ne 0, 1$, they satisfy the following three term recurrence formula:
\begin{equation} \label{ReR}
xM_{n}(x;\beta,c)=M_{n+1}(x;\beta,c)
+\frac{(c+1)n+\beta c}{1-c}M_{n}(x;\beta,c)
+\gamma_nM_{n-1}(x;\beta,c),
\end{equation} 
where 
\[
\gamma_n=\frac{cn(n+\beta-1)}{(1-c)^2},
\]
which can be explicitly in terms of hypergeometric series as\cite[(9.10.1)]{Koekoeketal} 
\begin{equation}
M_{n}(x;\beta,c)=\dfrac{c^n\,(\beta)_{n}}{(c-1)^n}
\,_{2}F_{1}\left(\begin{array}{c}-n,-x\\ \beta \end{array}; 1-\frac1c \right). \label{MRH}
\end{equation}

Next, we summarize some basic properties of Meixner orthogonal 
polynomials to be used in the sequel.

\begin{prop}
\label{S2-PROP-1} The following identities hold true for the 
Meixner polynomials: 

\begin{enumerate}
\item Second order difference equation. 
\begin{equation} \label{sodeme}
c(x+\beta)y_{n}(x+1)-(x(c+1)+\beta c)y_{n}(x)+x y_{n}(x-1)=n(c-1) y_{n}(x).
\end{equation} 
\item Structure relations. For every $n\in\mathbb{N}$, 
\begin{eqnarray}
\label{StruR1}(x+\beta)\Delta M_{n}(x;\beta,c)&=&
nM_{n}(x;\beta,c)+\frac{n( n+\beta-1)}{1-c} M_{n-1}(x;\beta,c),\\
\label{StruR2}x \nabla M_{n}(x;\beta,c)&=& nM_{n}(x;\beta,c)
+\frac{nc( n+\beta-1)}{1-c} M_{n-1}(x;\beta,c).
\end{eqnarray}

\item Squared norm. For every $n\in\mathbb{N}$, 
\begin{equation}
d_n^2=\| M_{n}(x;\beta,c)\|^{2}
=\left\langle {\bf u}^{\tt M},M^2_{n}(x;\beta,c)\right\rangle
=\frac{(\beta)_{n}c^{n}\, n!}{(1-c)^{\beta+2n}}. \label{Norm2}
\end{equation}
\item Value in the initial extreme of the orthogonality interval, 
\begin{equation}
M_{n}(0;\beta,c)=\dfrac{(\beta)_{n}c^n}{(c-1)^n}. 
\label{ViE}
\end{equation}
\item Forward and backward difference operators. For every $n, k\in\mathbb{N}$, 
\begin{equation}\label{FDO}
\Delta^{k}M_{n}(x;\beta,c)=(m-k+1)_k M_{n-k}(x;\beta+k,c),
\end{equation}
where $\Delta^{k} f(x)=\Delta^{k-1} \Delta f(x)$ for $k=1, 2,\dots$, and
$\nabla^{0} f(x)=f(x)$.
\item Mehler--Heine type formula \cite[eq. 35]{Dominici}
\begin{equation}
\lim_{n\rightarrow\infty }\frac{(1-c)^{n+\beta+z}M_{n}(z;\beta,c)}{(z-n+1)_n}=1,
\quad z\in \mathbb C\setminus \mathbb N.
\label{MHTF}
\end{equation}
\end{enumerate}
\end{prop}
To complete this section we present some useful results we need along the paper.
\begin{prop} (Christoffel-Darboux formula). 
Let $\left\{p_n\right\}_{n\in\mathbb N_0}$ be a sequence of monic polynomials 
orthogonal with respect to the linear functional $\bf u$. If we denote the 
$n$-th reproducing kernel by
\begin{equation} \label{Kernel1}
K_{n}(x,y):= \sum_{k=0}^{n-1}\frac{p_{k}(
x) p_{k}(y) }{\langle {\bf u}, p^2_k\rangle}.
\end{equation} 
Then, for all $n\in\mathbb{N}$, 
\begin{equation} \label{CDarb}
K_{n}(x,y) =\frac{1}{\langle {\bf u}, p^2_{n-1}\rangle}\frac{
p_{n}(x) p_{n-1}(y) -p_{n}(y) p_{n-1}(x) }{x-y}. 
\end{equation}
\end{prop}

Taking into account the inner product we have considered, then it 
is natural to consider the partial derivatives of $K_{n}(x,y) $ we
will use the following notation:
\begin{equation} \label{Kij}
{\mathscr K}_{n,1}^{(i,j) }(x,y) :=\sum_{k=0}^{n-1}
\frac{\nabla^{i}\, p_{k}(x)\nabla^{j}\, p_{k}(y) }{\langle {\bf u}, p^2_k\rangle},
\end{equation}
and
\[
{\mathscr K}_{n,2}^{(i,j) }(x,y) :=\sum_{k=0}^{n-1}
\frac{\Delta^{i}\, p_{k}(x)\Delta^{j}\, p_{k}(y) }{\langle {\bf u}, p^2_k\rangle}.
\]
\begin{prop} (cf. \cite[p. 454]{MR1106092})
The following identities hold:
\begin{equation} \label{Kernel0j-2}
{\mathscr K}_{n,1}^{(0,j)}(x,y)=\frac{j!}{d_{n-1}^{2}} 
\sum_{k=0}^{j} \left(\dfrac{M_{n}(x;\beta,c)\nabla^{k}M_{n-1}(y;\beta,c)
-M_{n-1}(x;\beta,c)\nabla^{k}M_{n}(y;\beta,c)}{k!\,(x-y+k)_{j+1-k}}\right), 
\end{equation} 

\begin{equation} \label{Kernel0j}
{\mathscr K}_{n,2}^{(0,j)}(x,y)=\frac{j!}{d_{n-1}^{2}} 
\sum_{k=0}^{j} \left(\dfrac{M_{n}(x;\beta,c)\Delta^{k}M_{n-1}(y;\beta,c)
-M_{n-1}(x;\beta,c)\Delta^{k}M_{n}(y;\beta,c)}{k!\,(x-y-j)_{j+1-k}}\right). 
\end{equation} 
\end{prop}

\begin{proof}
We are going to prove the first identity. 
After applying to \eqref{Kernel1} the difference operator $\nabla^j$ with
respect to $y$ we obtain 
\begin{equation}
{\mathscr K}_{n,1}^{(0,j)}(x,y)=\dfrac{1} {d_{n-1}^{2}}%
\Bigg(M_{n}(x;\beta,c)\nabla_y^j\left(\frac{M_{n-1}(x;\beta,c)}{x-y}\right)
-M_{n-1}(x;\beta,c)\nabla_y^j \left(%
\frac{M_{n}(y;\beta,c)}{x-y}\right)\Bigg). \label{K0jI}
\end{equation}
Using a analogue of the Leibnitz's rule
\begin{equation}
\nabla^n\left(f(x) g(x)\right)=\sum_{k=0}^{n}\binom{n}{k}
\nabla^k f(x)\nabla^{n-k} g(x-k) , \label{LR}
\end{equation}
and since for any positive integer $k$ we have
\begin{equation*}
\nabla_y^k\left(\frac{1}{x-y}\right)=\frac{k!}{(x-y)_{k+1}},
\end{equation*}
the result follows after a straightforward calculation. 
The proof of the second identity is analogous and it will be omitted.
Hence the result follows.
\end{proof}

\begin{prop}
The following identity holds for the kernel associated to the Meixner polynomials: 
\begin{equation*}
{\mathscr K}_{n,2}^{(j,j)}( 0,0)=\dfrac{j!(1-c)^{\beta+2j}}
{c^{j}(\beta)_{j}}\sum_{k=0}^{n-j-1}
\dfrac{( j+1)_{k}(\beta+j)_{k}}{(1)_{k}}\dfrac{c^{k}}{k!}.
\end{equation*}
\end{prop}

\begin{proof}
Having \eqref{ViE}-\eqref{Kernel1} into account and by definition\eqref{Kij} 
we have
\begin{eqnarray*}
{\mathscr K}_{n,2}^{(j,j)}( 0,0)&=&
\sum_{k=0}^{n-1} \frac{(k-j+1)^2_j(\beta+j)^2_j c^{2k-2j} (1-c)^{\beta+2k}}
{(c-1)^{2k-2j}(\beta)_j c^k k!}\\ &=&\frac{(1-c)^{2j+\beta}}{c^{2j}(\beta)^2_j}
\sum_{k=0}^{n-1} \frac{(k-j+1)^2_j(\beta)_k\, c^k}{k!}.
\end{eqnarray*}
By using some identities of the pochhammer symbols we obtain
\[\begin{split}
{\mathscr K}_{n,2}^{(j,j)}(0,0)=&
\frac{(1-c)^{2j+\beta}}{c^{2j} (\beta)^2_j}
\sum_{k=0}^{n-1}\frac{(k-j+1)^2_j(\beta)_kc^k}{(1)^2_{k-j}}\\[2mm]
=&\frac{(1-c)^{2j+\beta}}{c^{2j} (\beta)^2_j}
\sum_{k=0}^{n-1-j}\frac{(1)_{k+j}(\beta)_{k+j}c^{k+j}}{(1)^2_{k}},
\end{split}\]
from this expression it is a direct calculation to get the desired expression.
\end{proof}

\begin{cor}
\label{HyperGS} 
The following limit for the kernels associated to the Meixner polynomials holds:
\begin{equation*}
\lim_{n\rightarrow\infty }{\mathscr K}_{n,2}^{(j,j)}(0,0)
=\dfrac{j!(1-c)^{\beta+2j}}{c^{j}(\beta)_{j}} \,_{2}F_{1}\left(
\begin{array}{c}1+j, \beta+j\\1\end{array}; c \right).
\end{equation*}
\end{cor}
\section{The Sobolev-type Meixner polynomials}
We start this section introducing the Sobolev-type inner product \eqref{SobIP}
\begin{equation} \label{SobIP-1} 
\left\langle f,g\right\rangle_{\lambda,j,\ell}= 
\langle {\bf u}^{\tt M},fg\rangle
+\lambda \mathscr T^j f (\alpha)\mathscr T^{j}g(\alpha),
\end{equation} 
where ${\bf u}^{\tt M}$ is the Meixner linear operator, $j\in\mathbb N$, 
$\alpha\le 0$, and $\mathscr T$ is the operator $\nabla$ when $\ell=1$, and it 
is the operator $\Delta$ when $\ell=2$.

We denote by $\left\{{\mathscr M}^{j,\ell}_{n}(x;\beta,c;\lambda)\right\}_{n\in 
\mathbb N_0}$ 
the sequence of monic polynomials, orthogonal with respect to the 
inner product \eqref{SobIP-1}. These polynomials are said to be 
Sobolev-type Meixner polynomials.
\subsection{Connection formula and hypergeometric representation}
We first express the Sobolev-type Meixner polynomials 
in terms of the monic Meixner polynomials and the Kernel 
polynomials associated to the Meixner polynomials. 

Taking into account the Fourier expansion and using orthogonality 
conditions of $(M_n(x))$ and $({\mathscr M}^{j,\ell}_{n}(x;\beta,c;\lambda))$ 
we obtain (see cf. \cite[eq. (2.8)]{MR1106092})

\begin{equation}
{\mathscr M}^{j,\ell}_{n}(x;\beta,c;\lambda)=M_{n}(x;\beta,c)
-\dfrac{\lambda \mathscr{T}^{j}M_{n}(\alpha;\beta,c)}
{1+\lambda {\mathscr K}_{n,\ell}^{(j,j)}(\alpha ,\alpha)}
{\mathscr K}_{n,\ell}^{(0,j)}(x,\alpha), \ \ell=1,2. \label{FC}
\end{equation}
We can express the Sobolev-type Meixner polynomials in terms of the 
Meixner and their associated Kernel polynomials. 
Moreover, starting from \eqref{FC} and by using the recurrence 
relation of the Meixner polynomials \eqref{ReR} we have
\begin{equation}
{\mathscr M}^{j,\ell}_{n}(x;\beta,c;\lambda)={A}_{1,n}^{j,\ell}(x)
M_{n}(x;\beta,c)+{B}_{n,1}^{j,\ell}(x)M_{n-1}(x;\beta,c),
\label{SobP1}
\end{equation}
where 
\begin{align}
{A}_{1,n}^{j,1}(x)&=1-\dfrac{\lambda \mathscr{T}^{j}M_{n}(\alpha;\beta,c)}
{1+\lambda {\mathscr K}_{n,\ell}^{(j,j)}(\alpha ,\alpha)}
\frac{j!}{d_{n-1}^{2}} 
\sum_{k=0}^{j} \dfrac{\nabla^{k}M_{n-1}(y;\beta,c)}{k!\,(x-y+k)_{j+1-k}},\\
 {A}_{1,n}^{j,2}(x)&=1-\dfrac{\lambda \mathscr{T}^{j}M_{n}(\alpha;\beta,c)}
{1+\lambda {\mathscr K}_{n,\ell}^{(j,j)}(\alpha ,\alpha)}
 \frac{j!}{d_{n-1}^{2}} 
\sum_{k=0}^{j} \dfrac{\Delta^{k}M_{n-1}(y;\beta,c)}{k!\,(x-y-j)_{j+1-k}}, \\
{B}_{1,n}^{j,1}(x)&=\dfrac{\lambda \mathscr{T}^{j}M_{n}(\alpha;\beta,c)}
{1+\lambda {\mathscr K}_{n,\ell}^{(j,j)}(\alpha ,\alpha)}
\frac{j!}{d_{n-1}^{2}} 
\sum_{k=0}^{j} \frac{\nabla^{k}M_{n}(y;\beta,c)}{k!\,(x-y+k)_{j+1-k}},\\
{B}_{1,n}^{j,2}(x)&=\dfrac{\lambda \mathscr{T}^{j}M_{n}(\alpha;\beta,c)}
{1+\lambda {\mathscr K}_{n,\ell}^{(j,j)}(\alpha ,\alpha)}
\frac{j!}{d_{n-1}^{2}} 
\sum_{k=0}^{j} \frac{\Delta^{k}M_{n}(y;\beta,c)}{k!\,(x-y-j)_{j+1-k}}.
\end{align}

From these identities we can express the Sobolev-type Meixner 
polynomials in terms of hypergeometric series.
\begin{thm}
The monic Sobolev-type Meixner polynomial
${\mathscr M}^{j,\ell}_{n}(x;\beta,c;\lambda)$ has the following 
hypergeometric representation for $\ell=1,2$,
\begin{equation}
{\mathscr M}^{j,\ell}_{n}(x;\beta,c;\lambda)=\dfrac{(\beta)_{n-1}c^{n-1}}
{(c-1)^{n-1}} h_{n}^{\ell}(x)
\,_{3}F_{2}\left(\begin{array}{c}-n,-x, -f_{n}^{\ell}(x)\\ 
\beta, -f_{n}^{\ell}(x) -1 \end{array}; 1-\frac1c \right) , \label{MSPHR}
\end{equation}
where $f_{n}^{\ell}(x)$ is given in\eqref{fx} and
\begin{equation*}
h_{n}^{\ell}(x)=\frac{c (\beta+n-1)}{1-c} A_{1,n}^{j,\ell}(x)
-B_{1,n}^{j,\ell}(x).
\end{equation*}
\end{thm}

\begin{proof}
Taking into account $(-x)_{k}=0$ if $x<k$ 
as well as \eqref{MRH} and \eqref{SobP1} we deduce
\[\begin{split}
{\mathscr M}^{j,\ell}_{n}(x;\beta,c;\lambda)&=
\dfrac{c^n\,(\beta)_{n}}{(c-1)^n}\,A_{1,n}^{j,\ell}(x)
\sum_{k=0}^{n}\dfrac{(-n)_{k}(-x)_{k}}{(\beta)_{k}k!}
\left(1-\frac1c\right)^{k}\\
&+\dfrac{c^{n-1}\,(\beta)_{n-1}}{(c-1)^{n-1}}\,B_{1,n}^{j,\ell}(x)
\sum_{k=0}^{n-1}\dfrac{(1-n)_{k}(-x)_{k}}{(\beta)_{k}k!}
\left(1-\frac1c\right)^{k}.
\end{split}\]
By using the identity 
\begin{equation} \label{poid1} 
(a+k)(a)_k=a(a+1)_k,
\end{equation} 
we get
\begin{multline*}
{\mathscr M}^{j,\ell}_{n}(x;\beta,c;\lambda)=
\dfrac{c^n\,(\beta)_{n}}{(c-1)^n}\,A_{1,n}^{j,\ell}(x)
\sum_{k=0}^{n}\frac{(-n)_{k}(-x)_{k}}
{(\beta)_{k}k!}\left(1-\frac1c\right)^{k}\\
+\dfrac{c^{n-1}\,(\beta)_{n-1}}{(c-1)^{n-1}}\,\frac{B_{1,n}^{j,\ell}(x)}{n}
\sum_{k=0}^{n}
\dfrac{(n-k)(-n)_{k}(-x)_{k}}{(\beta)_{k}k!}
\left(1-\frac1c\right)^{k}.
\end{multline*}
Thus, we have
\begin{equation*}
{\mathscr M}^{j,\ell}_{n}(x;\beta,c;\lambda)\!=\!
\dfrac{(\beta)_{n-1}c^{n-1}}{(c-1)^{n-1}}
\frac{B_{1,n}^{j,\ell}(x)}{n}\!\!\sum_{k=0}^{n}
(f_{n}^{\ell}(x)-k+1)\frac{(-n)_{k}(-x)_{k}}
{(\beta)_{k}k!}\!\left(1-\frac1c\right)^{k},
\end{equation*}
where
\begin{equation}
f_{n}^{\ell}(x)=n-1-\frac{n c (\beta+n-1)A_{1,n}^{j,\ell}(x)}
{(1-c)B_{1,n}^{j,\ell}(x)}. \label{fx}
\end{equation}
And after a straightforward calculation and by using \eqref{poid1} 
with $a\to -f_{n}^{\ell}(x)-1$ the identity \eqref{MSPHR} follows. 
This completes the proof.
\end{proof}
\section{Second Order Linear Difference equation}
In this section we obtain a second order linear difference equation
that the sequence $\{{\mathscr M}^{j,\ell}_{n}(x;\beta,c;\lambda)\}_{n\geq 0}$ 
satisfies. 
In order to do that, we will find the \textit{ladder (creation and annihilation) 
operators}, using the connection formula \eqref{SobP1}, the three-term
recurrence relation\eqref{ReR}, and the structure relations \eqref{StruR1},
\eqref{StruR2} satisfied by them.

From \eqref{SobP1} and the recurrence relation\eqref{ReR} we 
deduce the following result
\begin{equation}\label{SobP2}
{\mathscr M}^{j,\ell}_{n-1}(x;\beta,c;\lambda)=
A_{2,n}^{j,\ell}(x) M_{n}(x;\beta,c)+B_{2,n}^{j,\ell}(x) 
M_{n-1}(x;\beta,c), 
\end{equation}
where
\[
A_{2,n}^{j,\ell}(x)=\frac{(c-1)B_{1,n-1}^{j,\ell}(x)}{(c+1)(n-1)+\beta c},
\quad \text{and}
\quad 
B_{2,n}^{j,\ell}(x)=A_{1,n-1}^{j,\ell}(x)+A_{2,n}^{j,\ell}(x)(1-x) .
\]
Applying the operator $\mathscr T$ to \eqref{SobP1}
and by using \eqref{LR} we have
\[\begin{split}
\mathscr T{\mathscr M}^{j,\ell}_{n}(x;\beta,c;\lambda&)=
M_{n}(x;\beta,c)\mathscr T A_{1,n}^{j,\ell}(x)
+A_{1,n}^{j,\ell}(x+(-1)^\ell)\mathscr TM_{n}(x;\beta,c)\\
+&M_{n-1}(x;\beta,c)\mathscr T
B_{1,n}^{j,\ell}(x)+A_{1,n}^{j,\ell}(x+(-1)^\ell) 
\mathscr TM_{n-1}(x;\beta,c).
\end{split}\]
Then, multiplying the previous expression by 
$x$ and using the structure relation\eqref{StruR2} if $\ell=1$ 
and $x+\beta$ and using the structure relation\eqref{StruR1} 
if $\ell=2$, and as well as the recurrence relation\eqref{ReR} 
we deduce the following expressions
\begin{equation}\label{DQ1}
x \nabla{\mathscr M}^{j,1}_{n}(x;\beta,c;\lambda)=C_{1,n}^{1}(x) 
M_{n}(x;\beta,c)+D_{1,n}^{1}(x) M_{n-1}(x;\beta,c), 
\end{equation}
\begin{equation}\label{DQ1-2}
(x+\beta)\Delta {\mathscr M}^{j,2}_{n}(x;\beta,c;\lambda)=C_{1,n}^{2}(x) 
M_{n}(x;\beta,c)+D_{1,n}^{2}(x) M_{n-1}(x;\beta,c), 
\end{equation}
\begin{equation}\label{DQ2-1}
x \nabla{\mathscr M}^{j,1}_{n-1}(x;\beta,c;\lambda)=
C_{2,n}^{1}(x)M_{n}(x;\beta,c)+D_{2,n}^{1}(x)M_{n-1}(x;\beta,c), 
\end{equation}
and
\begin{equation}\label{DQ2}
(x+\beta)\Delta{\mathscr M}^{j,2}_{n-1}(x;\beta,c;\lambda)=
C_{2,n}^{2}(x)M_{n}(x;\beta,c)+D_{2,n}^{2}(x)M_{n-1}(x;\beta,c), 
\end{equation}
respectively, where all the coefficients can be computed explicitly.
Moreover, from \eqref{SobP1}-\eqref{SobP2} for $\ell=1, 2$ we have
\begin{equation*}
\Theta _{n}(x;\ell)M_{n}(x;\beta,c)=B_{2,n}^{\ell}(x) 
{\mathscr M}^{j,\ell}_{n}(x;\beta,c;\lambda) -B_{1,n}^{\ell}(x) 
{\mathscr M}^{j,\ell}_{n-1}(x;\beta,c;\lambda),
\end{equation*}
and
\begin{equation*}
\Theta _{n}(x;\ell) M_{n-1}(x;\beta,c)=A_{1,n}^{\ell}(x)
{\mathscr M}^{j,\ell}_{n-1}(x;\beta,c;\lambda) -A_{2,n}^{\ell}(x)
{\mathscr M}^{j,\ell}_{n}(x;\beta,c;\lambda), 
\end{equation*}
where
\begin{equation*}
\Theta _{n}(x;\ell)=\det 
\begin{pmatrix}
A_{1,n}^{\ell}(x)& B_{1,n}^{\ell}(x)\\[3mm] 
A_{2,n}^{\ell}(x)& B_{2,n}^{\ell}(x)
\end{pmatrix},\qquad \ell=1,2.
\end{equation*}
After replacing the above in\eqref{DQ1}, \eqref{DQ1-2} 
and \eqref{DQ2-1}, \eqref{DQ2}, we conclude
\[
\left(\tilde{\Theta}_{n}(x;\ell)\mathscr{T}+\Lambda _{2,n}^{(1)}(x;\ell)
\right)\left[{\mathscr M}^{j,\ell}_{n}(x;\beta,c;\lambda)\right]
=\Lambda _{1,n}^{(1)}(x;\ell) {\mathscr M}^{j,\ell}_{n-1}(x;\beta,c;\lambda).
\]
and
\begin{equation*}
\left(\tilde{\Theta}_{n}(x;\ell)\mathscr{T} 
+\Lambda _{1,n}^{(2)}(x;\ell)\right)\left[{\mathscr M}^{j,\ell}_{n-1}(x;\beta,c;\lambda)\right]
=\Lambda _{2,n}^{(2)}(x;\ell){\mathscr M}^{j,\ell}_{n}(x;\beta,c;\lambda),
\end{equation*}
respectively, where
\begin{equation}\label{Coef1}
\tilde{\Theta}_{n}(x;\ell)=\left\{
\begin{array}{r@{\ \text{if}\ } l} 
x \Theta _{n}(x;\ell),&\ell=1,\\[3mm]
(x+\beta)\Theta _{n}(x;\ell),&\ell=2,\end{array}\right.
\end{equation}
and
\begin{equation} \label{Coef2}
\Lambda _{j,n}^{(k)}(x;\ell)=(-1)^{k}\det \begin{pmatrix}
C_{k,n}^{\ell}(x)& A_{\nu,n}^{j,\ell}(x)\\[3mm] 
D_{k,n}^{\ell}(x)& B_{\nu,n}^{j,\ell}(x)
\end{pmatrix}, \ \ \nu=1,2,\ \ k=1,2,\quad \ell=1,2. 
\end{equation}

\begin{prop}
Let $\left\{{\mathscr M}^{j,\ell}_{n}(x;\beta,c;\lambda)
\right\}_{n\in\mathbb N_0}$  be the sequence of monic 
Sobolev-type Meixner polynomials defined by \eqref{MSPHR} 
and let $\mathscr I$ be the identity operator. Then, the ladder
(destruction and creation) operators $\mathfrak{a}$, 
$\mathfrak{a}^{\dagger}$ are defined by
\begin{equation*}
\mathfrak{a}=\tilde{\Theta}_{n}(x;\ell)\mathscr T+\Lambda
_{2,n}^{(1)}(x;\ell)\mathscr I,
\end{equation*}
\begin{equation*}
\mathfrak{a}^{\dagger }=\tilde{\Theta}_{n}(x;\ell)\mathscr T
+\Lambda _{1,n}^{( 2)}(x;\ell)\mathscr I,
\end{equation*}
which verify%
\begin{equation}\label{DO}
\mathfrak{a}({\mathscr M}^{j,\ell}_{n}(x;\beta,c;\lambda))=
\Lambda_{1,n}^{(1)}(x;\ell) {\mathscr M}^{j,\ell}_{n-1}(x;\beta,c;\lambda), 
\end{equation}
\begin{equation*}
\mathfrak{a}^{\dagger }({\mathscr M}^{j,\ell}_{n-1}(x;\beta,c;\lambda))
=\Lambda_{2,n}^{( 2)}(x;\ell) {\mathscr M}^{j,\ell}_{n}(x;\beta,c;\lambda) ,
\end{equation*}
where $\mathscr I$ is the identity operator, 
$\tilde{\Theta}_{n}(x;\ell) $ and $\Lambda _{j,n}^{(k)}(x;\ell)$ 
with $j, k,\ell=1, 2$ are given in\eqref{Coef1}-\eqref{Coef2}.
\end{prop}

\begin{thm}
The monic Sobolev-type Meixner polynomials sequence, 
which is orthogonal with respect to the inner product \eqref{SobIP-1}, fulfills 
the second order difference equation:
\begin{equation}
\mathcal{F}_{n}(x;\ell)\mathscr{T}^{2}y(x)
+\mathcal{G}_{n}(x;\ell)\mathscr{T}y(x)
+\mathcal{H}_{n}(x;\ell)y(x)=0, \label{HEq}
\end{equation}
where
\begin{equation*}
\mathcal{F}_{n}(x;\ell)=\frac{\tilde{\Theta}_{n}(x;\ell) 
\tilde{\Theta}_{n}(x+(-1)^\ell;\ell)}{\Lambda_{1,n}^{(1)}( x+(-1)^\ell;\ell)},
\end{equation*}
\[\begin{split}
\mathcal{G}_{n}(x;\ell)=&\frac{\tilde{\Theta}_{n}(x;\ell)}{\Lambda_{1,n}^{(1)}(x+(-1)^\ell;\ell)}
\left(\mathscr{T}\tilde{\Theta}_{n}(x;\ell)-\frac{\tilde{\Theta}_{n}(x;\ell)
\mathscr{T}\Lambda _{1,n}^{(1)}(x;\ell)}{\Lambda_{1,n}^{(1)}(x;\ell)} \right.
\\[2mm]&+\Lambda _{2,n}^{(1)}(x+(-1)^\ell;\ell)\Bigg)+\frac{\tilde{\Theta}_{n}(x;\ell)\Lambda
_{1,n}^{( 2)}(x;\ell)}{\Lambda_{1,n}^{(1)}(x;\ell)},
\end{split}\]
and
\[\begin{split}
\mathcal{H}_{n}(x;\ell)=&\frac{\tilde{\Theta}_{n}(x;\ell)\mathscr{T}
\Lambda _{2,n}^{(1)}(x;\ell)}{\Lambda_{1,n}^{(1)}( x+(-1)^\ell;\ell)}
-\dfrac{\tilde{\Theta}_{n}(x;\ell)\Lambda _{2,n}^{(1)}(x;\ell) 
\mathscr{T}\Lambda _{1,n}^{(1)}(x;\ell)}{\Lambda_{1,n}^{(1)}(x;\ell) 
\Lambda _{1,n}^{(1)}(x+(-1)^\ell;i)}\\ &+\dfrac{\Lambda_{1,n}^{( 2)}(x;\ell) 
\Lambda_{2,n}^{(1)}(x;\ell)}{\Lambda_{1,n}^{(1)}(x;\ell)}
-\Lambda _{2,n}^{( 2)}(x;\ell),
\end{split}\]
where $\tilde{\Theta}_{n}(x;\ell) $ and 
$\Lambda _{j,n}^{(k)}(x;\ell) $ with $j,k,\ell=1,2$ are given in 
\eqref{Coef1} and \eqref{Coef2}.
\end{thm}

\begin{proof}
From \eqref{DO} we have
\begin{equation*}
\frac{1}{\Lambda_{1,n}^{(1)}(x;\ell)}\mathfrak{a}
({\mathscr M}^{j,\ell}_{n}(x;\beta,c;\lambda))=
{\mathscr M}^{j,\ell}_{n-1}(x;\beta,c;\lambda).
\end{equation*}
Next, applying the operator $\mathfrak{a}^{\dagger }$ to both members 
of the previous expression, we get
\begin{equation*}
\mathfrak{a}^{\dagger}\left[ \frac{1}{\Lambda_{1,n}^{(1)
}(x;\ell)}\mathfrak{a}({\mathscr M}^{j,\ell}_{n}(x;\beta,c;\lambda)
)\right] =\Lambda _{2,n}^{( 2)}(x;\ell)
{\mathscr M}^{j,\ell}_{n}(x;\beta,c;\lambda) .
\end{equation*}
Thus, by using the definitions of the operators $\mathfrak{a}$ and 
$\mathfrak{a}^{\dagger}$, taking into account the identity 
\begin{equation*}
\mathscr{T}\left\{\dfrac{f(x)}{g(x)}\right\} =
\frac{g(x)\mathscr{T}f(x)-f(x)\mathscr{T}g(x)}
{g(x) g(x+(-1)^\ell)},
\end{equation*}
and after tedious calculations we obtain\eqref{HEq}. Hence the result follows.
\end{proof}
\section{The $(2j+3)$-term recurrence relation}
In this section we find the $(2j+3)$-term recurrence relation
that the sequence of monic Sobolev-type Meixner 
polynomials \eqref{MSPHR} fulfill. For this purpose, 
we use the fact, which is a straightforward consequence of 
\eqref{SobIP}, that the multiplication operator by $\left\langle
x-\alpha \right\rangle_{i}^{j+1}$ is a symmetric operator with respect to
such a discrete Sobolev inner product. Indeed, for any $p, q\in\mathbb{P}$
we have for $\ell=1$ 
\begin{equation} \label{POk}\begin{split}
\left\langle (x-\alpha)_{j+1}p(x),q(x)\right\rangle_{\lambda,j,\ell}
=&\left\langle (x-\alpha)_{j+1}{\bf u}^{M}, p(x)q(x)\right\rangle \\[2mm]
=&\left\langle p(x), (x-\alpha)_{j+1}q(x)\right\rangle_{\lambda,j,\ell},
\end{split}\end{equation}
and for $\ell=2$ 
\begin{equation} \label{POk-1}\begin{split}
\left\langle (x-\alpha-j)_{j+1}p(x),q(x)\right\rangle_{\lambda,j,\ell}
=&\left\langle (x-\alpha-j)_{j+1}{\bf u}^{M}, p(x)q(x)\right\rangle\\[2mm]
=&\left\langle p(x), (x-\alpha-j)_{j+1}q(x)\right\rangle_{\lambda,j,\ell}.
\end{split}\end{equation}
Taking these identities into account and by using the three-term recurrence
relation\eqref{ReR} we can state the following result.
\begin{lem}
The following identities related to the monic Sobolev-type Meixner polynomials 
hold:
\begin{eqnarray}\label{XaC}
\hspace{-5mm}(x-\alpha)_{j+1}{\mathscr M}^{j,1}_{n}(x;\beta,c;\lambda)=&
\mathcal{A}^{1}_{n}(x) M_{n}(x;\beta,c)+\mathcal{B}^{1}_{n}(x)  M_{n-1}(x;\beta,c), \\[2mm]
\label{XaC-1}\hspace{-5mm}(x-\alpha-j)_{j+1}
{\mathscr M}^{j,2}_{n}(x;\beta,c;\lambda)=&\mathcal{A}^{2}_{n}(x) 
M_{n}(x;\beta,c)+\mathcal{B}^{2}_{n}(x)  M_{n-1}(x;\beta,c), 
\end{eqnarray}
where $\mathcal{A}^{\ell}_{n}(x)$, $\mathcal{B}^{\ell}_{n}(x)$
are polynomials which can be computed explicitly.
\end{lem}

\begin{thm}
Let $\lambda \in\mathbb{R}_{+}$,  and $j\in\mathbb{N}$,  
let $\left\{{\mathscr M}^{j,\ell}_{n}(x;\beta,c;\lambda)\right\}_{n\in\mathbb N_0}$ 
be the sequence of monic Sobolev-type Meixner polynomials defined 
by \eqref{MSPHR}.  

Then, the norm of these polynomials fulfills the following identity:
\begin{equation} \label{NormSP}
\|{\mathscr M}^{j,\ell}_{n}(x;\beta,c;\lambda)\|_{\lambda,j,\ell}^{2}
=\|M_{n}(x;\beta,c)\|^{2}+b_{n}^{j,\ell}\|M_{n-1}(x;\beta,c)\|^{2}, 
\end{equation}
where 
\begin{equation} \label{bnk}
b_{n}^{j,\ell}=\frac{\lambda}{\|M_{n-1}(x;\beta,c)\|^{2}}
\frac{\mathscr{T}^{j}M_{n}(\alpha;\beta,c)\mathscr{T}^{j}M_{n}(\alpha;\beta,c)}
{1+\lambda {\mathscr K}_{n,\ell}^{(j,j)}(\alpha ,\alpha)}\ge 0.
\end{equation}
\end{thm}
\begin{proof}
We will consider the $\ell=1$ case. The $\ell=2$ case is analogous.

By the property of orthogonality of Sobolev-type Meixner polynomials we have
\begin{equation*}
\|{\mathscr M}^{j,\ell}_{n}(x;\beta,c;\lambda)\|_{\lambda,j,\ell}^{2}=
\left\langle{\mathscr M}^{j,\ell}_{n}(x;\beta,c;\lambda),(x-\alpha)_{j+1}\pi_{n-j-1}(x)
\right\rangle_{\lambda,j,\ell},
\end{equation*}
for any monic polynomial $\pi$ of degree $n-j-1$. 
From \eqref{POk} we have
\[
\left\langle {\mathscr M}^{j,\ell}_{n}(x;\beta,c;\lambda), 
(x-\alpha)_{j+1}\pi(x\right\rangle_{\lambda,j,\ell}=
\left\langle (x-\alpha)_{j+1}{\mathscr M}^{j,\ell}_{n}(x;\beta,c;\lambda), \pi(x)
\right\rangle_{\lambda,j,\ell}.
\]
By using the connection formula \eqref{XaC} and taking
into account that $\mathcal{A}^{\ell}_{n}(x) $ is a monic polynomial
of degree  $j+1$ and $\mathcal{B}^{\ell}_{n}(x) $ is a
polynomial of degree  $j$ with the leading coefficient $b_{n}^{j,\ell}$ 
we deduce
\[\begin{split}
\|{\mathscr M}^{j,\ell}_{n}(x;\beta,c;\lambda)\|_{\lambda,j,\ell}^{2} =&\left\langle 
(x-\alpha)_{j+1}{\mathscr M}^{j,\ell}_{n}(x;\beta,c;\lambda), \pi(x)
\right\rangle_{\lambda,j,\ell}\\
=&\left\langle {\bf u}^{\tt M}, \!\mathcal{A}^{\ell}_{n}(x) M_{n}(x;\beta,c)\pi(x)
\right\rangle\!+\!\left\langle {\bf u}^{\tt M},\! \mathcal{B}^{\ell}_{n}(x) 
M_{n-1}(x;\beta,c) \pi(x)\right\rangle \\ 
=&\left\langle {\bf u}^{\tt M}, M_{n}(x;\beta,c)\, x^{n}\right\rangle
+b_{n}^{j,\ell}\left\langle {\bf u}^{\tt M}, M_{n-1}(x;\beta,c)\,x^{n-1}\right\rangle,
\end{split}\]
which coincides with \eqref{NormSP}.
\end{proof}
\begin{rem} 
Observe that a direct consequence is 
\begin{equation} \label{relnorm} 
\frac{\|{\mathscr M}^{j,\ell}_{n}(x;\beta,c;\lambda)\|_{\lambda,j,\ell}^{2}}
{\|M_{n}(x;\beta,c)\|^{2}}=
\frac{1+\lambda {\mathscr K}_{n+1,\ell}^{(j,j)}(\alpha ,\alpha)}
{1+\lambda {\mathscr K}_{n,\ell}^{(j,j)}(\alpha ,\alpha)}
.
\end{equation} 
\end{rem} 
\begin{thm}[$(2j+3)$-term recurrence relation]
Let $\lambda \in\mathbb{R}_{+}$, $j\in\mathbb{N}_0$. Then, the monic 
Sobolev-type Meixner orthogonal polynomials sequence with respect to 
the inner product \eqref{SobIP-1} satisfies the following $(2j+3)$-term 
recurrence relation:
\begin{equation*}
(x-\alpha)_{j+1}{\mathscr M}^{j,1}_{n}(x;\beta,c;\lambda)=
{\mathscr M}^{j,1}_{n+j+1}(x;\beta,c;\lambda)
+\sum_{k=n-j-1}^{n+j}c_{n,k}^{j,1}\, {\mathscr M}^{j,1}_{k}(x;\beta,c;\lambda),
\end{equation*}
and 
\begin{equation*}
(x-\alpha-j)_{j+1}{\mathscr M}^{j,2}_{n}(x;\beta,c;\lambda)=
{\mathscr M}^{j,2}_{n+j+1}(x;\beta,c;\lambda)
+\hspace{-2mm}\sum_{k=n-j-1}^{n+j}\hspace{-2mm} c_{n,k}^{j,1}\, {\mathscr 
M}^{j,2}_{k}(x;\beta,c;\lambda),
\end{equation*}
where the constant coefficients $c_{n,k}^{j,\ell}$ can be explicitly computed  
for $\ell=1, 2$.
\end{thm}

\begin{proof}
In such a case we will consider the $\ell=2$ case. 

Since the Sobolev-type Meixner polynomials form a basis in 
$L_2(\langle \cdot,\cdot \rangle_{\lambda,j,\ell})$, 
if we consider the Fourier expansion of 
$(x-\alpha-j)_{j+1}{\mathscr M}^{j,\ell}_{n}(x;\beta,c;\lambda)$ 
in terms of the Sobolev-type Meixner polynomials, then
\begin{equation*}
(x-\alpha-j)_{j+1}{\mathscr M}^{j,\ell}_{n}(x;\beta,c;\lambda)
={\mathscr M}^{j,\ell}_{n+j+1}(x;\beta,c;\lambda)
+\sum_{k=0}^{n+j}c^{j,\ell}_{n,k} {\mathscr M}^{j,\ell}_{k}(x;\beta,c;\lambda),
\end{equation*}
Thus, by using the property of orthogonality of the sequence 
$\{{\mathscr M}^{j,\ell}_{n}(x;\beta,c;\lambda)\}$ we obtain
\begin{equation*}
c^{j,\ell}_{n,k}=\frac{\left\langle(x-\alpha-j)_{j+1}{\mathscr M}^{j,\ell}_{n}(x;\beta,c;\lambda), 
{\mathscr M}^{j,\ell}_{k}(x;\beta,c;\lambda)\right\rangle_{\lambda,j,\ell}}
{\|{\mathscr M}^{j,\ell}_{k}(x;\beta,c;\lambda)\|_{\lambda,j,\ell}^{2}},\ \ k=0,... ,n+j.
\end{equation*}
Using \eqref{POk-1} and the property of orthogonality 
of  $\{{\mathscr M}^{j,\ell}_{n}(x;\beta,c;\lambda)\}$, we deduce 
that $c^{j,\ell}_{n,k}=0$ for $k=0, \ldots, n-j-2$. 
Observe that the rest of the coefficients can be computed by using again the 
same orthogonality conditions. The proof of the another identity is similar and 
it will be omited.
\end{proof}

\section{Mehler-Heine type formula}

The main result of this section will be to establish Mehler--Heine 
type formula of the polynomial ${\mathscr M}^{j,\ell}_{n}(x;\beta,c;\lambda)$ 
for the  $\alpha\le 0$ case.  Let us see the following result.

\begin{lem} \label{lem:5}
Let $\beta, c\in\mathbb C$, with $|c|<1$ and $-\beta\not \in \mathbb N$, 
and let $m$ be a positive integer.
Then, the following limit holds:
\begin{equation}
\lim_{n\rightarrow\infty }\frac{(\beta)_{n}n^m c^{n}}{(n-1)!}=0, \label{Lm1}
\end{equation}
\end{lem}

\begin{proof}
If we use the identity \cite[p. 23]{DRain} 
\begin{equation*}
\Gamma(z)=\lim_{n\rightarrow\infty }\frac{(n-1)!n^{z}}{( z)_{n}},
\end{equation*}
we deduce
\begin{equation} \label{Lim2}
\lim_{n\rightarrow\infty }\frac{(\beta)_{n}n^m c^{n}}{(n-1)!}
=\frac{1}{\Gamma(\beta)}\lim_{n\rightarrow\infty}
n^{\beta+m}c^{n}. 
\end{equation} 
Therefore, if Re$(\beta+m)>0$ then applying to \eqref{Lim2} L'H\^opital's 
rule several times we obtain the desired result; otherwise the limit is zero. 
Hence the result holds.
\end{proof}

\begin{lem}
Let $\beta, c\in\mathbb R$, with $|c|<1$ and $\beta$ is not a negative integer, 
and let $k, j$ be integers, with $0\leq k\leq j$. If we set $\alpha=0$ in\eqref{FC}. 
Then, the following limits hold:
\begin{equation}
\lim_{n\rightarrow\infty }{A}_{1,n}^{j,\ell}(x)=1\quad 
\mbox{and}\quad \lim_{n\rightarrow\infty }{B}_{1,n}^{j,\ell}(x)=0,\qquad 
\ell=1,2. \label{Limanbn}
\end{equation}
\end{lem}

\begin{proof}
By starting with \eqref{FC} and using \eqref{Kernel0j} we obtain
\[
{A}_{i,n}^{j,\ell}(x)\!=1-\dfrac{\lambda\mathscr{T}^{j}M_{n}(\alpha;\beta,c)}
{1+\lambda {\mathscr K}_{n,\ell}^{(j,j)}(\alpha,\alpha)}
\frac{j!}{d_{n-1}^{2}} \sum_{k=0}^{j} \frac{\mathscr{T}^{k}M_{n-1}(\alpha;\beta,c)}
{k!\,(x-\alpha+k)_{j+1-k}},
\]
and
\[
{B}_{i,n}^{j,\ell}(x)\!=\lambda \dfrac{\mathscr{T}^{j}M_{n}(\alpha;\beta,c)}
{1+\lambda {\mathscr K}_{n,\ell}^{(j,j)}(\alpha,\alpha)}\frac{j!}{d_{n-1}^{2}} 
\sum_{k=0}^{j} \frac{\mathscr{T}^{k}M_{n}(\alpha;\beta,c)}
{k!\,(x-\alpha+k)_{j+1-k}}.
\]
where $i=1, 2$. Then, to prove this result it is enough to check
\begin{equation*}
\lim_{n\rightarrow\infty }(1-c)^n a^{j,\ell}_{k,n}=\lim_{n\rightarrow
\infty }(1-c)^n b^{j,\ell}_{k,n}=0,\ \ k=0, 1,\dots, n,\ \ i=1, 2, \quad \ell=1,2,
\end{equation*}
where
\begin{equation*}
a^{j,\ell}_{k,n}=\dfrac{\lambda\, j!\, \mathscr{T}^{j}M_{n}(\alpha;\beta,c)
\mathscr{T}^{k}M_{n-1}(\alpha;\beta,c)}
{k!(1+\lambda {\mathscr K}_{n,\ell}^{(j,j)}(\alpha,\alpha)) d_{n-1}^{2}},
\end{equation*}
and%
\begin{equation*}
b^{j,\ell}_{k,n}=\dfrac{\lambda\, j!\, \mathscr{T}^{j}M_{n}(\alpha;\beta,c)
\mathscr{T}^{k}M_{n}(\alpha;\beta,c)} {k!(1+\lambda 
{\mathscr K}_{n,\ell}^{(j,j)}(\alpha,\alpha)) d_{n-1}^{2}}.
\end{equation*}
After a straightforward calculation, by using \eqref{MHTF} we have that for any 
$0\le k\le j$
\[
\mathscr{T}^{k}M_{n}(\alpha;\beta,c)\approx 
\mathscr{T}^{j}M_{n}(\alpha;\beta,c),
\]
for $n$ large, and since $c-1<c$, then by using Lemma \ref{lem:5} it is clear that 
both limits related to such coefficients tend to zero. Hence  we deduce 
\eqref{Limanbn}.
\end{proof}

\begin{thm}
Let $\beta, c\in\mathbb R$, with $0<c<1$ and $\beta$ is not a negative integer,
and let $m$ be a positive integer.
Then, we have
\begin{equation}
\lim_{n\rightarrow\infty }\frac{(1-c)^{n+\beta+z}{\mathscr M}^{j,\ell}_{n}(z;\beta,c;\lambda)}
{(z-n+1)_n}=1, \quad z\in \mathbb C\setminus \mathbb N.
\label{MHQ1n}
\end{equation}
uniformly on compact subsets of the complex plane.
\end{thm}

\begin{proof}
Multiplying \eqref{SobP1} by the factor 
$(1-c)^{n+\beta+z}/(z-n+1)_n$  we have
\[\begin{split} 
\frac{(1-c)^{n+\beta+z}{\mathscr M}^{j,\ell}_{n}(z;\beta,c;\lambda)}{(z-n+1)_n}
=&A_{i,n}^{(j,\ell)}(z)\frac{(1-c)^{n+\beta+z}M_{n}(z;\beta,c)}{(z-n+1)_n}\\
& +B_{1,n}^{(j,\ell)}(z)\frac{(1-c)^{n+\beta+z}M_{n-1}(z;\beta,c)}{(z-n+1)_n}.
\end{split}\]
Then, applying the previous Lemma as well as the \eqref{MHTF} 
we arrived to the desired result.
\end{proof}

\begin{rem}
Observe that we can extend some of the previous results even for 
$c\in \mathbb C$ so that $|c-1|<|c|<1$, or even into a wider region of 
the complex plane taking into account \eqref{meorex}.
\end{rem}

Finally, we show some graphical experiments of the limit function in 
\eqref{MHQ1n} for several values of $n$ using Mathematica software at 
the masspoint $\alpha=0$,  see Figures \ref{GMH50}, \ref{GMH70}, and 
\ref{GMH100}.

\begin{figure}[!hbt]
\begin{center}
\begin{tikzpicture}[domain=-1:3.5, samples=100, scale=1.4]
\draw[-stealth](-2,0) -- (5,0);
\draw[-stealth] (0,-2.5) --(0,2.5);
\foreach \x/\xtext in {-1.5/-1,1.5/1,3/2,4.5/3}
\draw[shift={(\x,0)}] (0pt,2pt) -- (0pt,-2pt) node[below] {\tiny $\xtext$};
\foreach \y/\ytext in {-1/-5,-2/-10,1/5,2/10}
\draw[shift={(0,\y)}] (2pt,0pt) -- (0pt,0pt) node[left] {\tiny $\ytext$};	
\draw[thin, blue] plot[domain=-1.5:1.5] (\x,{(-4.76837*(2*\x/3)
+1.68835*(2*\x/3)^2+3.62293*(2*\x/3)^3+0.557285*(2*\x/3)^4
-0.756343*(2*\x/3)^5-0.377349*(2*\x/3)^6)/5}) 
plot[domain=1.5:3] (\x,{(2.81982+2.612*(-1.5+(2*\x/3))-12.014147*(-1.5+(2*\x/3))^2-
11.9867755*(-1.5+(2*\x/3))^3+2.45458276*(-1.5+(2*\x/3))^4+6.27385*(-1.5+(2*\x/3))^5)/5})
plot[domain=3:4.75] (\x,{(-8.81194-11.6873*(-2.5+(2*\x/3))+34.2792*(-2.5+(2*\x/3))^2+
 52.4764*(-2.5+(2*\x/3))^3+7.3129*(-2.5+(2*\x/3))^4-22.674*(-2.5+(2*\x/3))^5-
 14.2863*(-2.5+(2*\x/3))^6)/5});
\draw[thick, red, densely dashed] plot[domain=-1.5:1.5] 
(\x,{(-4.94613*(2*\x/3)+1.53952*(2*\x/3)^2
+3.79685*(2*\x/3)^3+0.749453*(2*\x/3)^4- 0.73173*(2*\x/3)^5-0.418988*(2*\x/3)^6)/5})
plot[domain=1.5:3] (\x, {(3.16606+3.13296*(-1.5+(2*\x/3))-13.2745*(-1.5+(2*\x/3))^2-
 14.2823*(-1.5+(2*\x/3))^3+1.78213*(-1.5+(2*\x/3))^4+7.08437*(-1.5+(2*\x/3))^5)/5})
 plot[domain=3:4.7] (\x,{(-10.6191-14.9172*(-2.5+(2*\x/3))+40.0872*(-2.5+(2*\x/3))^2+
 66.3134*(-2.5+(2*\x/3))^3+14.2151*(-2.5+(2*\x/3))^4-25.8969*(-2.5+(2*\x/3))^5-
 19.1988*(-2.5+(2*\x/3))^6)/5});
\end{tikzpicture}
\end{center}
\caption[]{Limit function in \eqref{MHQ1n} for $n=50$, (red color) 
left member and right member in blue. \\
Data: $\beta = 7$, $c = 1/5$, $\lambda = 10^{-21}$ and $j=2$.}
\label{GMH50}
\end{figure}
\begin{figure}[!hbt]
\begin{center}
\begin{tikzpicture}[domain=-1:3.5, samples=100, scale=1.4]
\draw[-stealth](-2,0) -- (5,0);
\draw[-stealth] (0,-2.5) --(0,2.5);
\foreach \x/\xtext in {-1.5/-1,1.5/1,3/2,4.5/3}
\draw[shift={(\x,0)}] (0pt,2pt) -- (0pt,-2pt) node[below] {\tiny $\xtext$};
\foreach \y/\ytext in {-1/-5,-2/-10,1/5,2/10}
\draw[shift={(0,\y)}] (2pt,0pt) -- (0pt,0pt) node[left] {\tiny $\ytext$};	
\draw[thin, blue] plot[domain=-1.5:1.5] (\x,{(-4.76837*(2*\x/3)
+1.68835*(2*\x/3)^2+3.62293*(2*\x/3)^3+0.557285*(2*\x/3)^4
-0.756343*(2*\x/3)^5-0.377349*(2*\x/3)^6)/5})
plot[domain=1.5:3] (\x,{(2.81982+2.612*(-1.5+(2*\x/3))
-12.0141*(-1.5+(2*\x/3))^2-11.9868*(-1.5+(2*\x/3))^3 
+2.45458*(-1.5+(2*\x/3))^4+6.27385*(-1.5+(2*\x/3))^5+
 2.06209*(-1.5+(2*\x/3))^6)/5})
plot[domain=3:4.7] (\x,{(-8.81194-11.6873*(-2.5+(2*\x/3))
+34.2792*(-2.5+(2*\x/3))^2+52.4764*(-2.5+(2*\x/3))^3+7.3129*(-2.5+(2*\x/3))^4
-22.674*(-2.5+(2*\x/3))^5-14.2863*(-2.5+(2*\x/3))^6)/5});
\draw[thick, red, densely dashed] plot[domain=-1.5:1.5] 
(\x, {(-4.89299*(2*\x/3)+1.58583*(2*\x/3)^2
+3.74654*(2*\x/3)^3+0.690363*(2*\x/3)^4-0.741205*(2*\x/3)^5-0.406968*(2*\x/3)^6)/5})
plot[domain=1.5:3] (\x,{(3.05651+2.96358*(-1.5+(2*\x/3))-12.883*(-1.5+(2*\x/3))^2-
 13.54*(-1.5+(2*\x/3))^3+2.02636*(-1.5+(2*\x/3))^4+6.83922*(-1.5+(2*\x/3))^5+
 2.54017*(-1.5+(2*\x/3))^6)/5})
 plot[domain=3:4.67] (\x,{(-10.0216 - 13.8214*(-2.5+(2*\x/3))+38.2187*(-2.5+(2*\x/3))^2+
 61.6529*(-2.5+(2*\x/3))^3+11.7133*(-2.5+(2*\x/3))^4 - 24.9542*(-2.5+(2*\x/3))^5 - 
 17.5339*(-2.5+(2*\x/3))^6)/5});
\end{tikzpicture}
\caption[]{Limit function in \eqref{MHQ1n} for $n=70$, (red color) 
left member and right member in blue.  \\
Data: $\beta = 7$, $c = 1/5$, $\lambda = 10^{-21}$ and $j=2$.}
\label{GMH70}
\end{center}
\end{figure} 
\begin{figure}[!hbt]
\begin{center}
\begin{tikzpicture}[domain=-1:3.5, samples=100, scale=1.4]
\draw[->](-2,0) -- (5,0);
\draw[->] (0,-2.5) --(0,2.5);
\foreach \x/\xtext in {-1.5/-1,1.5/1,3/2,4.5/3}
\draw[shift={(\x,0)}] (0pt,2pt) -- (0pt,-2pt) node[below] {\tiny $\xtext$};
\foreach \y/\ytext in {-1/-5,-2/-10,1/5,2/10}
\draw[shift={(0,\y)}] (2pt,0pt) -- (0pt,0pt) node[left] {\tiny $\ytext$};	
\draw[thin, blue] plot[domain=-1.5:1.5] (\x,{(-4.76837*(2*\x/3)
+1.68835*(2*\x/3)^2+3.62293*(2*\x/3)^3+0.557285*(2*\x/3)^4
-0.756343*(2*\x/3)^5-0.377349*(2*\x/3)^6)/5})
plot[domain=1.5:3] (\x,{(2.81982+2.612*(-1.5+2*\x/3))
-12.0141*(-1.5+(2*\x/3))^2-11.9868*(-1.5+(2*\x/3))^3 
+2.45458*(-1.5+(2*\x/3))^4+6.27385*(-1.5+(2*\x/3))^5+
 2.06209*(-1.5+(2*\x/3))^6)/5})
plot[domain=3:4.7] (\x,{(-8.81194-11.6873*(-2.5+(2*\x/3))
+34.2792*(-2.5+(2*\x/3))^2+52.4764*(-2.5+(2*\x/3))^3
+7.3129*(-2.5+(2*\x/3))^4-22.674*(-2.5+(2*\x/3))^5
-14.2863*(-2.5+(2*\x/3))^6)/5});
\draw[thick, red, densely dashed] plot[domain=-1.5:1.5] 
(\x,{(-0.3979869*(2*\x/3)^6-0.7469248*(2*\x/3)^5
+0.6483362*(2*\x/3)^4+3.709099*(2*\x/3)^3+1.618465*(2*\x/3)^2-4.854422*(2*\x/3))/5})
plot[domain=1.5:3] (\x,{(2.38156*(-1.5+(2*\x/3))^6+6.66221*(-1.5+(2*\x/3))^5+2.17871
  *(-1.5+(2*\x/3))^4-13.0327*(-1.5+(2*\x/3))^3-12.6068
  *(-1.5+(2*\x/3))^2+2.84833*(-1.5+(2*\x/3))+2.98037)/5})
plot[domain=3:4.67] (\x,{(-16.4372 *(-2.5+(2*\x/3))^6-24.2505
   *(-2.5+(2*\x/3))^5+10.1574 *(-2.5+(2*\x/3))^4+58.5725
   *(-2.5+(2*\x/3))^3+36.9361 *(-2.5+(2*\x/3))^2-13.1017
   *(-2.5+(2*\x/3))-9.62071)/5});
\end{tikzpicture}
\caption[]{Limit function in \eqref{MHQ1n} for $n=100$, (red color) 
left member and  right member in blue. \\
Data: $\beta=7$, $c = 1/5$, $\lambda = 10^{-21}$ and $j=2$.}
\label{GMH100}
\end{center}
\end{figure}
\vspace{6pt} 




\section*{Acknowledgments:} 

The research of RSCS was funded by Agencia Estatal de Investigaci\'on of Spain, grant num- ber PGC-2018-096504-B-C33. 
The work of ASL is supported by Direcci\'on General de Investigaci\'on e Innovaci\'on, Consejer\'ia de Educaci\'on e Investigaci\'on 
of the Comunidad de Madrid (Spain) and Universidad de Alcal\'a under grant CM/JIN/2021-014, Proyectos de I+D para J\'ovenes 
Investigadores de la Universidad de Alcal\'a 2021 .
We are grateful for the exhaustive work of the referees. Their comments and suggestions have improved the presentation of the manuscript.
\end{document}